

\baselineskip=14pt
\parskip=10pt

\font\eightrm=cmr8 
\font\eighttt=cmtt8
\magnification=\magstephalf

\def\1{{\overline{1}}}
\def\2{{\overline{2}}}
\parindent=0pt
\overfullrule=0in

\def\frac#1#2{{#1 \over #2}}
\bf
\centerline
{
Using Noonan-Zeilberger Functional Equations to enumerate (in Polynomial Time!) 
}
\centerline
{
Generalized Wilf classes
}
\rm
\bigskip
\centerline{ {\it
Brian NAKAMURA${}^1$ and
Doron 
ZEILBERGER}\footnote{$^1$}
{\eightrm  \raggedright
Department of Mathematics, Rutgers University (New Brunswick),
Hill Center-Busch Campus, 110 Frelinghuysen Rd., Piscataway,
NJ 08854-8019, USA.
{\eighttt [bnaka, zeilberg]  at math dot rutgers dot edu} ,
\hfill \break
{\eighttt http://www.math.rutgers.edu/\~{}bnaka/} \quad , \quad {\eighttt http://www.math.rutgers.edu/\~{}zeilberg/} .
\hfill\break
First Written: Sept. 7, 2012.  
Supported in part by the USA National Science Foundation.
}
}

{\bf Very Important} 

This article is accompanied by the Maple packages

$\bullet$
{\eighttt http://www.math.rutgers.edu/\~{}zeilberg/tokhniot/P123} \quad ,

$\bullet$
{\eighttt http://www.math.rutgers.edu/\~{}zeilberg/tokhniot/F123} \quad ,

$\bullet$
{\eighttt http://www.math.rutgers.edu/\~{}zeilberg/tokhniot/P1234} \quad , 

$\bullet$
{\eighttt http://www.math.rutgers.edu/\~{}zeilberg/tokhniot/F1234} \quad ,

$\bullet$
{\eighttt http://www.math.rutgers.edu/\~{}zeilberg/tokhniot/P12345} \quad , 

$\bullet$
{\eighttt http://www.math.rutgers.edu/\~{}zeilberg/tokhniot/F12345} \quad ,

$\bullet$
{\eighttt http://www.math.rutgers.edu/\~{}zeilberg/tokhniot/P123456} \quad , 

to be described below.  Lots of sample input and output files can be seen at:

{\tt http://www.math.rutgers.edu/\~{}zeilberg/mamarim/mamarimhtml/Gwilf.html} \quad \quad .

{\bf Introduction}

Recall  that the {\it reduction} of a finite list of $k$, say, distinct (real) numbers $[a_1, a_2 , \dots, a_k]$
is the unique permutation  $\sigma=[\sigma_1, \dots, \sigma_k]$,
of $\{1, \dots , k\}$ 
such that $a_1$ is the $\sigma_1$-th largest element in the list,
$a_2$ is the $\sigma_2$-th largest element in the list, etc. In other words
$[a_1, a_2 , \dots, a_k]$ and $\sigma$ are ``order-isomorphic''. For example,
the reduction of $[6,3,8,2]$ is $[3,2,4,1]$ and the reduction of $[\pi,\gamma,e,\phi]$ is
$[4,1,3,2]$.

Given a permutation $\pi=\pi_1 \dots \pi_n$ and another permutation $\sigma=[\sigma_1 , \dots , \sigma_k]$
(called a {\it pattern}), we denote by $N_\sigma(\pi)$ the number of instances $1 \leq i_1 < \dots < i_k \leq n$ such that the reduction of
$\pi_{i_1} \dots \pi_{i_k}$ is $\sigma$.

For example, if $\pi=51324$ then 

$N_{[1,2,3]}(\pi)=2$ (because $\pi_2 \pi_3 \pi_5=134$ and $\pi_2 \pi_4 \pi_5=124$ reduce to $[1,2,3]$).

$N_{[1,3,2]}(\pi)=1$ (because $\pi_2 \pi_3 \pi_4=132$  reduces to $[1,3,2]$).

$N_{[2,1,3]}(\pi)=1$ (because $\pi_3 \pi_4 \pi_5=324$  reduces to $[2,1,3]$).

$N_{[2,3,1]}(\pi)=0$ (because none of the $10$ length-three subsequences of $\pi$ reduces to $231$).

$N_{[3,1,2]}(\pi)=5$ (because $\pi_1 \pi_2 \pi_3=513$ 
and  $\pi_1 \pi_2 \pi_4=512$
and  $\pi_1 \pi_2 \pi_5=514$
and  $\pi_1 \pi_3 \pi_5=534$
and  $\pi_1 \pi_4 \pi_5=524$ all reduce to $[3,1,2]$) .

$N_{[3,2,1]}(\pi)=1$ (because $\pi_1 \pi_3 \pi_4=532$  reduces to $[3,2,1]$).

Of course the sum of $N_\sigma(\pi)$ over all $k$-permutations $\sigma$ is ${{n} \choose {k}}$.

Fixing a pattern $\sigma$, the set of permutations $\pi$ for which $N_\sigma(\pi)=0$ (we say that $\pi$ {\it avoids} $\sigma$)
is called the {\it Wilf class} of $\sigma$, and more generally, given a set of patterns
$S$, the set of permutations for which $N_\sigma(\pi)=0$ for all $\sigma \in S$, is the Wilf class
of that set.  The first {\it systematic} study  of {\it enumerating} Wilf classes
was undertaken in the pioneering paper by Rodica Simion and Frank Schmidt [SiSc].

The general question is extremely difficult (see [Wiki] and [Bo3]) and ``explicit''
answers are only known for few short patterns (and sets of patterns), the increasing patterns
$[1, 2, \dots , k]$, and a few other {\it West-equivalent} to them, giving the same enumeration.
For example, even for the pattern $[1,3,2,4]$ ({\tt http://oeis.org/A061552}) the best
known algorithm takes exponential time in $n$, and it is very possible that that's the
best that one can do.

But for those patterns $\sigma$ for which we know how to enumerate their Wilf classes, most importantly the
increasing patterns $[1, \dots , k]$, it makes sense to ask the more general question:

Given a pattern $\sigma$, and a positive integer $r$, find a ``formula'', or at least a polynomial-time
algorithm (thus {\it answering} the question in the sense of Herb Wilf[Wil])
that inputs a positive integer $n$ and outputs the number of permutations $\pi$ of $\{1, \dots, n\}$
for which $N_\sigma(\pi)=r$. We call such a class a {\it generalized} Wilf class.

Ideally, we would like to have, given  a pattern $\sigma$,  an explicit formula, in $n$ and $q$, for
the generating function ($S_n$ denotes the set of permutations of $\{1, \dots, n\}$)
$$
A_\sigma(q,n):=\sum_{\pi \in S_n} q^{N_\sigma(\pi)} \quad,
$$
then, for any fixed $r$, the sequence of coefficients of $q^r$ in  $A_\sigma(q,n)$ would give
the sequence enumerating permutations with {\it exactly} $r$ occurrences of the pattern $\sigma$.

In fact, for patterns of length $\leq 2$ there are nice answers. Trivially
$$
A_{[1]}(q,n):=n!\, q^n \quad ,
$$
and almost-trivially (or at least classically)
$$
A_{[2,1]}(q,n):=(1)\,(1+q) \, \dots \, (1+q+ \dots +q^{n-1}) \, = \, [n]! \quad,
$$
the famous ``$q$-analog'' of $n!$. 
But things start to get complicated for patterns of length $3$.

\vfill\eject

{\bf Past Work}

For a very lucid and extremely engaging introduction to the subject,
as well as the state-of-the-art in 2004, we strongly recommend
Mikl\'os B\'ona's {\it masterpiece} [Bo3].

In [NZ], John Noonan and the second-named author initiated 
a {\it functional equations}-based approach for
enumerating generalized Wilf classes.
In order to illustrate it, they
reproved John Noonan's[N] combinatorially-proved result
that the  number of permutations of length $n$ with exactly
{\it one} occurrence of the pattern $[1,2,3]$ equals
$\frac{3}{n} { {2n} \choose {n-3}}$.
Recently a {\it proof from the book} of this result was
given by Alexander Burstein[Bu] (see also [Z1]).

In [NZ] it was conjectured that the  number of permutations
of length $n$ with
{\it two} occurrences of the pattern $[1,2,3]$ equals
$\frac{59n^2+117n+100}{2n(2n-1)(n+5)} { {2n} \choose {n-4}}$.
This conjecture was proved by Markus Fulmek[F],
using Dyck paths.

In [NZ] it was also  conjectured that the  number of permutations
of length $n$ with
{\it one} occurrence of the pattern $[1,3,2]$ equals
$\frac{n-2}{2n} { {2n-2} \choose {n-1}}$.
This conjecture was proved by Mikl\'os B\'ona[Bo1],
who later proceeded to prove[Bo2] the interesting fact
that the sequences enumerating permutations 
with exactly $r$ occurrences of $[1,3,2]$ 
is  $P$-recursive (i.e. satisfies a homogeneous linear recurrence
with polynomial coefficients) for {\it every} $r$. In fact he
proved the stronger result that the generating functions are
always algebraic. This was vindicated by Toufik Mansour and
Alek Vainshtein[MV] who gave an efficient algorithm to actually
compute these generating functions, and they used it to
find explicit expressions for $1 \leq r \leq 5$.

Another interesting but different ``functional equation'' approach,
for patterns of length three, was developed by Firro and Mansour[FM].

{\bf This Project}

But so far, practically nothing is known for patterns of length
larger than three and $r>0$. In this paper we will modify
the approach of [NZ] in order to generate, in {\it polynomial  time},
such sequences for increasing patterns of {\it any} length
$[1, \dots ,k]$. That method can be extended to the patterns
$[1, \dots ,k-2,k,k-1]$ and possibly other families, but here
we will only discuss increasing patterns.

Using the new algorithm to compute sufficiently many
terms,
we were able to conjecture explicit formulas, in $n$, for the number of
permutations of length $n$ with exactly $r$ occurrences of the pattern $[1,2,3]$,
for $5 \leq r \leq 7$, extending Fulmek's[F] conjectures for $r=3$ and $r=4$. 
We believe that the {\it enumeration schemes}, that our algorithms generate,
should enable our computers to conjecture {\it holonomic} representations
for the more general quantities (see below), that once guessed, should be
amenable to automatic rigorous proving  in the {\it holonomic paradigm}[Z2],
using Christoph Koutschan's[K] far-reaching extensions and powerful implementations.
But since these conjectures are {\it certainly} true, and their formal proof 
would (probably) not yield any new insight, we don't think that it is
worth the trouble to actually carry out the gory details, wasting
both humans' time (it would require quite a bit of daunting programming) and
the computers' time (it would take a very long time, due to the complexity of
the schemes).

Now let's recall the {\it Noonan-Zeilberger Functional Equation Approach}.

{\bf The Noonan-Zeilberger Functional Equation Approach}

The {\it starting point} of the
Noonan-Zeilberger[NZ] approach
for enumerating generalized Wilf classes is to
derive  a {\it functional equation}. 
Let's review it
with the simplest non-trivial case, that of the length-$3$ 
increasing pattern $[1,2,3]$.

In addition to the variable $q$, introduce $n$ extra {\it catalytic variables} $x_1, \dots, x_n$, and
define the {\it weight} of a permutation 
$\pi=\pi_1 \dots \pi_n$ of length $n$ by
$$
weight(\pi):=q^{N_{[1,2,3]}(\pi)} \prod_{i=1}^{n} 
x_i^{|\{1\leq a <b \leq n; \pi_a=i<\pi_b\}|} \quad ,
$$
(as usual, for any set $A$, $|A|$ denotes the number of elements of $A$).
For example,
$$
weight(12345)=q^{10} x_1^4 x_2^3 x_3^2 x_4 \quad ,
$$
$$
weight(54321)=1 \quad ,
$$
$$
weight(21354)=q^{4} x_2^3 x_1^3 x_3^2= q^{4} x_1^3 x_2^3  x_3^2 \quad .
$$

Let's define the polynomial in the $n+1$ variables
$$
P_n(q; x_1, \dots, x_n):=
\sum_{\pi \in S_n} weight(\pi) \quad .
$$
Let $\pi=\pi_1 \dots  \pi_n$ be a typical permutation of length $n$. Suppose $\pi_1=i$.
Note that the number of occurrences of the pattern $[1,2,3]$ in $\pi$ equals
the number of occurrences of that pattern in the beheaded permutation $\pi_2 \dots \pi_n$
{\bf plus} the number of the patterns $[1,2]$ in the beheaded permutation $\pi_2 \dots \pi_n$
where the ``1'' is $i+1$, or $i+2$, or $\dots$ or $n$.  
Let $\pi'$ be the reduction to $\{1, \dots , n-1\}$ of 
that beheaded permutation.
We see that
$$
weight(\pi)=x_i^{n-i} weight(\pi')
\,\, 
\vert_{\, x_i \rightarrow qx_{i+1}
\,\, , \,\,  x_{i+1} \rightarrow qx_{i+2}, \,\, \dots \,\, ,  \,\,
x_{n-1} \rightarrow qx_{n}}  \quad .
$$
The factor of $x_i^{n-i}$ is because converting $\pi'$ from a permutation of $\{1, \dots, n-1\}$ to
a permutation of $\{1, \dots, i-1, i+1, \dots ,n \}$, and sticking an $i$ at the front
introduces $n-i$ new $[1,2]$ patterns where the ``$1$'' is $i$.
This gives the {\it Noonan-Zeilberger Functional Equation} 
$$
P_n(q; x_1 , \dots , x_n)=
\sum_{i=1}^{n} x_i^{n-i} P_{n-1}(q; x_1, \dots, x_{i-1}, qx_{i+1}, \dots, qx_n) \quad .
\eqno(NZFE1)
$$

Having found $P_n(q;x_1, \dots, x_n)$, we set the ``catalytic'' variables $x_1, \dots, x_n$ all to $1$
and get
$$
f_n(q):=A_{[1,2,3]}(q,n)=P_n(q; 1, 1, \dots, 1 ) \quad.
$$
Even though this is an ``exponential-time'' (and memory!) algorithm, it is much faster than
the direct weighted counting of all the $n!$ permutations, and we were able to
explicitly compute them through $n=20$. 

{\eightrm This is implemented in  procedure {\eighttt fn(n,q)} in {\eighttt P123}. Procedure
{\eighttt L20(q);} gives the pre-computed sequence of {\eighttt fn(n,q)} for n between 1 and 20 . }

Here are the first few terms:
$$
f_1(q)=1 \quad, \quad
f_2(q)=2 \quad, \quad
f_3(q)=q+5 \quad, \quad
f_4(q)={q}^{4}+3\,{q}^{2}+6\,q+14 \quad, \quad
$$

$$
f_5(q)={q}^{10}+4\,{q}^{7}+6\,{q}^{5}+9\,{q}^{4}+7\,{q}^{3}+24\,{q}^{2}+27\,q+42 \quad ,
$$

$$
f_6(q)={q}^{20}+5\,{q}^{16}+8\,{q}^{13}+6\,{q}^{12}+6\,{q}^{11}+16\,{q}^{10}+12\,{q}^{9}+24\,{q}^{8}+32\,{q}^{7}+37\,{q}^{6}+54\,{q}^{5}+74\,{q}^{4}+70\,{q}^{3}
+133\,{q}^{2}+110\,q+132 \quad ,
$$

$$
f_7(q)=
{q}^{35}+6\,{q}^{30}+10\,{q}^{26}+10\,{q}^{25}+8\,{q}^{23}+13\,{q}^{22}+30\,{q}^{21}+10\,{q}^{20}+32\,{q}^{19}+18\,{q}^{18}+62\,{q}^{17}+74\,{q}^{16}+24
\,{q}^{15}+100\,{q}^{14}
$$
$$
+130\,{q}^{13}+104\,{q}^{12}+162\,{q}^{11}+191\,{q}^{10}+232\,{q}^{9}+260\,{q}^{8}+320\,{q}^{7}+387\,{q}^{6}+395\,{q}^{5}+507\,{q
}^{4}+461\,{q}^{3}+635\,{q}^{2}+429\,q+429 \quad ,
$$

$$
f_8(q)=
{q}^{56}+7\,{q}^{50}+12\,{q}^{45}+15\,{q}^{44}+10\,{q}^{41}+16\,{q}^{40}+40\,{q}^{39}+18\,{q}^{38}+47\,{q}^{36}+38\,{q}^{35}+68\,{q}^{34}+60\,{q}^{33}
$$
$$
+58\,{q}^{32}+66\,{q}^{31}+154\,{q}^{30}+138\,{q}^{29}+115\,{q}^{28}+156\,{q}^{27}+252\,{q}^{26}+324\,{q}^{25}+228\,{q}^{24}+288\,{q}^{23}+537\,{q}^{22}
$$
$$
+466\,{q}^{21}+546\,{q}^{20}+656\,{q}^{19}+682\,{q}^{18}+1004\,{q}^{17}+1047\,{q}^{16}+886\,{q}^{15}+1494\,{q}^{14}+1456\,{q}^{13}+1580\,{q}^{12}
$$
$$
+1818\,{q}^{
11}+2077\,{q}^{10}+2182\,{q}^{9}+2389\,{q}^{8}+2544\,{q}^{7}+2864\,{q}^{6}+2570\,{q}^{5}+3008\,{q}^{4}+2528\,{q}^{3}+2807\,{q}^{2}+1638\,q+1430
\quad .
$$

For $f_n(q)$ for $9 \leq n \leq 20$ see:

{\tt http://www.math.rutgers.edu/\~{}zeilberg/tokhniot/oP123d} . 

Using this data, the computer easily finds {\it rigorously-proved}
explicit expressions for the first six moments (about the mean) of the random variable ``number of occurrences of the
pattern $[1,2,3]$'', and from them verifies that, at least up to the sixth moment, this random variable is
{\it asymptocally normal}, as humanly proved (for {\it all} patterns) by Mikl\'os B\'ona[Bo4].  See:

{\tt http://www.math.rutgers.edu/\~{}zeilberg/tokhniot/oP123a} . 

{\bf The ``Perturbation'' Approach}

The equations of {\it quantum field theory} are (usually) impossible to solve exactly, but
physicists got around it by devising clever ``approximate'' methods using
{\it perturbation} expansions, that only use the first few terms in a potentially ``infinite''
(and intractable) series, but that suffice for all practical purposes, using {\it Feynman diagrams}.

Of course, we are {\it enumerators}, and we want {\it exact} results, but suppose we only
want to know the sequences enumerating permutations with exactly $s$ occurrences of the pattern $[1,2,3]$
for $s \leq r$ for some relatively small $r$, rather than for $r={{n} \choose {3}}$, provided by the
full 
$$
f_n(q)=P_n(q\, ; \, 1 \, [n \,\,\, times \, ]) \quad .
$$

In the original article [NZ], for $r=0$, Noonan and Zeilberger  simply plugged-in $q=0$ and
$x_1 = \dots =x_n=1$, getting a simple enumeration scheme, that proved, for the $n$-th time,
the classical result that the number of permutations of length $n$ that avoid the pattern $123$
equals the Catalan number $(2n)!/(n!(n+1)!)$. For $r=1$, they differentiated Eq. $(NZFE1)$ with  respect
to $q$, using the multivariable calculus {\it chain rule}, and then plugged-in $q=0$ and 
$x_1 = \dots =x_n=1$. For $r=2$ they did it again, but
this turned out to be, for larger $r$,  a {\it Rube Goldberg} nightmare, {\it even for a computer}.

Here is a much  easier way!

Recall that you are really only interested in $f_n(q)=P_n(q \, ; \, 1 \, [ n \,\,\, times] \,)$.
Plugging it into (NZFE1) gives
$$
P_n(q\,; \, 1 \, [ n \,\,\, times \,]) \, = \,
\sum_{i=1}^{n}  P_{n-1}(q \, ; \, 1 \, [ \, i-1 \,\,\, times \, ] \, , \, q \, [\, n-i \,\,\, times \,]) \quad .
$$
This forces us to put-up with expressions of the form
$$
P_{a_0+a_1}(q; \, 1\, [\,a_0 \,\,\, times] \, , \, q \,[ \, a_1 \,\,\, times]) \quad .
$$
Plugging this into $(NZFE1)$ yields
$$
P_{a_0+a_1}(q; \, 1 \, [\, a_0 \,\,\, times] \, , \, q \,[ \, a_1 \,\,\, times])=
\sum_{i=1}^{a_0}  P_{a_0+a_1-1}(q; \, 1 \, [ \, i-1 \,\,\, times] \, , \,q \,[ \, a_0-i \,\,\, times] \, , \,q^2 [a_1 \,\,\, times]) 
$$
$$
+\sum_{i=1}^{a_1}  
q^{a_1-i} P_{a_0+a_1-1}(q;\, 1 \, [ \, a_0 \,\,\, times] \, , \, q \, [ \, i-1 \,\,\, times] \, , \, q^2 \, [ \, a_1-i \,\,\, times])
\quad .
$$
This forces us, in turn, to consider expressions of the form
$$
P_{a_0+a_1+a_2}(q; \, 1 \, [ \, a_0 \,\,\, times]\, , \, q \, [ \, a_1 \,\,\, times] \, , \, q^2 \, [ \, a_2 \,\,\, times] ), \quad
$$
that would force us to further consider expressions of the form
$$
P_{a_0+a_1+a_2+a_3}(q; \, 1 \, [\, a_0 \,\,\, times \, ] \, , \, q \, [ \, a_1 \,\,\, times] \, ,  \,  
\, q^2 \, [ \, a_2 \,\,\, times] \, , \,q^3 \, [ \, a_3 \,\,\, times] ) \quad ,
$$
etc. etc., leading to an {\bf exponential explosion} in both time and memory.

But, if we are only interested in the first $r$ 
coefficients of $f_n(q)$, then we can take advantage of the
{\it trivial} (you prove it!) but {\bf crucial} lemma.

{\bf Crucial Lemma}: For $s >r+1$, the coefficients of $q^0,q^1, \dots, q^r$ of
$$
P_{a_0+a_1+ \dots +a_s}(q;\, 1 \, [ \, a_0 \,\,\, times \,] 
\, , \, \dots q^{s-1} \, , \, [\, a_{s-1} \,\,\, times \,]
\, , \,  q^s \, [\, a_s \,\,\, times \,] )
$$
$$
- \, P_{a_0+a_1+ \dots +a_s}(q;\, 1 \, [\, a_0 \,\,\, times] \, , \, \dots \, , \, 
 q^{r} \, [a_r \,\,\, times \,] \, , \,
q^{r+1} \, [\, a_{r+1}+a_{r+2}+ \dots + a_{s} \,\,\, times \,] )
$$
all vanish. In other words, if  for any polynomial 
$p(q)$ in $q$, $p^{(r)}(q)$ denotes the polynomial of
degree $r$ obtained by  ignoring all powers of $q$ larger than $r$,
and letting $CHOP_r[\, p(q) \, ]:=p^{(r)}(q)$,
$(NZFE1)$ becomes 
(below, let  $n:=a_1+ \dots +a_r+a_{r+1}$, and for any expression $R$ and positive integer $k$, $R\$k$,
denotes $R \dots R [ k \,\,\, times]$ for example $q^2\$3$ means $q^2,q^2,q^2$) 
$$
P^{(r)}_{n}(q;\, 1\$a_0 \, , \, 
q\$a_1\, , \, 
\dots  \, , \, q^{r}\$a_{r}, q^{r+1}\$a_{r+1} )
$$
$$
=CHOP_r (
\sum_{i=1}^{a_0} P^{(r)}_{n-1}(q;\, 1\$(i-1)\, , \, q\$(a_0-i) \, ,  \,
q^{2}\$a_{1} \, , \, \dots,  \, 
\,
q^{r}\$a_{r-1} \, , \,
q^{r+1}\$(a_r+a_{r+1}) )
$$
$$
+
\sum_{i=1}^{a_1} q^{a_1-i+a_2+ \dots + a_{r+1}}
P^{(r)}_{n-1}(q; \, 1\$a_0
\, , \, q\$(i-1)
\, , \,
q^{2}\$(a_{1}-i) \, ,  \,
\dots \, ,  \, 
q^{r}\$a_{r-1} \, ,  \, 
q^{r+1}\$(a_r+a_{r+1}))
$$
$$
+
\sum_{i=1}^{a_2} q^{2(a_2-i+a_3+ \dots + a_{r+1})}
P^{(r)}_{n-1}(q; \, 1\$a_0  \, , \, q\$a_1
\, , \,
q^{2}\$(i-1),
q^{3}\$(a_{2}-i) \, ,  \,
\dots \, , \,
q^{r}\$a_{r-1}
\, ,  \, q^{r+1}\$(a_r+a_{r+1}))
$$
$$
+\, \dots  \dots \dots
$$
$$
+
\sum_{i=1}^{a_{r+1}} q^{(r+1)(a_{r+1}-i)} 
P^{(r)}_{n-1}(q; \, 1 \$a_0  \, , \, 
q\$a_1  \, ,  
\, q^{2}\$a_2  \, \dots \, ,
q^{r}\$a_{r} \, ,  \, 
 \, q^{r+1}\$(a_{r+1}-1) ) ) \quad .
$$

Now note that, because of the $CHOP_r$ operator in front, many terms
automatically disappear, because of the powers of $q$ in front.
The {\it bottom line} is that the computer can automatically generate a 
scheme for computing
the degree-$r$ polynomials in $q$,
$$
F_r(a_0, \dots, a_{r+1})(q) \, := \,
P^{(r)}_{a_0+ \dots + a_{r+1}} (q;\, 1 \, [\, a_0 \,\,\, times \, ] \, , \, q \, [ \, a_1 \,\,\, times] \, , \, \dots \, , \,
q^{r+1} \, [ \, a_{r+1} \,\,\, times \,]) \quad ,
$$
with $a_0+ \dots + a_{r+1}=n$ and $a_0, \dots , a_{r+1} \geq 0$.
The number of such quantities is the coefficient of $z^n$ in $1/(1-z)^{r+2}$ that equals
$(-1)^{r+2}{ {-(r+2)} \choose {n}}={ {r+n+1} \choose {r+1}}$ terms.
So each iteration involves $O(n^{r+1})$  evaluations and hence $O(n^{r+2})$
additions and doing it $n$ times yields an $O(n^{r+3})$ algorithm for finding our object of desire,
the degree $r$ polynomial in $q$:
$$
f^{(r)}_n(q)\, = \, F_r(n \, , \, 0 \, [\, r+1 \,\,\, times \, ])(q) \quad .
$$
Having found the scheme, the very same computer (or a different one),
may use it to generate as many terms as desired.

{\bf The Maple package P123}

The Maple package P123 downloadable from

{\eighttt http://www.math.rutgers.edu/\~{}zeilberg/tokhniot/P123}  \quad ,

implements the functional equation $(NZFE1)$ and easily  generated the first $25$ terms of the 
enumerating sequences for $0 \leq r \leq 7$.
With this data, it empirically verified the
already-known results for the number of permutations with exactly $r$ occurrences of the
pattern $[1,2,3]$ for $0 \leq r \leq 2$, and made conjectures for $r \leq 7$ as follows.
Let $a_r(n)$ be the number of permutations of length $n$ with exactly $r$ occurrences
of the pattern $[1,2,3]$.

$$
a_0(n)=
2\,{\frac { \left( 2\,n-1 \right) !}{ \left( n-1 \right) !\, \left( n+1 \right) !}} \quad .
$$

$$
a_1(n)=
6\,{\frac { \left( 2\,n-1 \right) !}{ \left( n-3 \right) !\, \left( n+3 \right) !}} \quad .
$$

$$
a_2(n)={\frac { \left( 2\,n-2 \right) !}{ \left( n-4 \right) !\, \left( n+5 \right) !}} \,\cdot \,(59\,{n}^{2}+ 117\,n + 100 ) \quad .
$$

$$
a_3(n)=
{\frac { \left( 2\,n-3 \right) !}{ \left( n-5 \right) !\, \left( n+7 \right) !}} \, \cdot
4\,n \left(113\,{n}^{3} + 506\,{n}^{2} +937\,n + 1804\right) \quad .
$$

$$
a_4(n)=\frac{(2n-4)!}{(n-4)!(n+9)!} \, \cdot
$$
$$
\left ( 3561\,{n}^{8}+3126\,{n}^{7}-46806\,{n}^{6}+12384\,{n}^{5}-659091\,{n}^{4}+2630634\,{n}^{3}+5520576\,{n}^{2}+26283456\,n-39191040 \right )
\quad .
$$

$$
a_5(n)=\frac{(2n-5)!}{(n-5)!(n+11)!} \, \cdot 
$$
$$
( \, 
26246\,{n}^{10}+136646\,{n}^{9}-115872\,{n}^{8}+22524\,{n}^{7}-9648450\,{n}^{6}+71304534\,{n}^{5}
$$
$$
+381205612\,{n}^{4}+1607633896\,{n}^{3}+2800103664\,{n}^
{2}+3611692800\,n-32891443200  ) \quad .
$$

$$
a_6(n)=\frac{(2n-6)!}{(n-6)!(n+13)!} \cdot
$$
$$
( \, 193311\,{n}^{12}+2349954\,{n}^{11}+13035003\,{n}^{10}+95151030\,{n}^{9}+406430793\,{n}^{8}+2889552582\,{n}^{7}
$$
$$
+14335663329\,{n}^{6}+60005854890\,{n}^{5}+
313010684796\,{n}^{4}+1025692693464\,{n}^{3}
$$
$$
+1283595375168\,{n}^{2}-6909513045120\,n-28177269120000 \, ) \quad .
$$

$$
a_7(n)=\frac{(2n-7)!}{(n-5)!(n+15)!} \cdot
$$
$$
(\,
1386032\,{n}^{16}+13111080\,{n}^{15}+22526480\,{n}^{14}+355187760\,{n}^{13}-1654450096\,{n}^{12}+10534951680\,{n}^{11}
$$
$$
+15797223760\,{n}^{10}-305671694640\,{n}^{9}+3750695521216\,{n}^{8}-26631101348520\,{n}^{7}
$$
$$
-86395090065440\,{n}^{6}-636425872408320\,{n}^{5}+3647384624274048\,{n}^{4}
$$
$$
+11386434230674560\,{n}^{3}+103032675524966400\,{n}^{2}-157858417817856000\,n-763734137886720000 \,)
\quad .
$$

{\bf Enumerating Permutations with r occurrences of the pattern [1,2,3,4] for small r
via a Noonan-Zeilberger Functional Equation}

In addition to the variable $q$, we now introduce $2n$ extra {\it catalytic variables} 
$x_1, \dots, x_n$, and $y_1, \dots, y_n$, and
define the {\it weight} of a permutation 
$\pi=\pi_1 \dots \pi_n$ of length $n$ by
$$
weight(\pi):=q^{N_{[1,2,3,4]}(\pi)} \prod_{i=1}^{n} 
x_i^{|\{1\leq a <b <c\leq n \, ; \,  \pi_a=i<\pi_b<\pi_c\}|}  \cdot
y_i^{|\{1\leq a <b \leq n \, ; \, \pi_a=i<\pi_b\}|}  \quad .
$$
For example,
$$
weight([1,2,3,4,5,6])=q^{15} x_1^{10} x_2^6 x_3^3 x_4 y_1^5 y_2^4 y_3^3 y_4^2 y_5 \quad ,
$$
$$
weight([6,5,4,3,2,1])=1 \quad ,
$$
$$
weight([3,4,5,6,1])=q x_3^3 x_4 y_3^3 y_4^2 y_5 \quad .
$$

Let's define the polynomial in the $2n+1$ variables
$$
P_n(q; x_1, \dots, x_n; y_1, \dots, y_n):=
\sum_{\pi \in S_n} weight(\pi) \quad .
$$
Let $\pi=\pi_1  \dots , \pi_n$ be a typical permutation of length $n$. Suppose $\pi_1=i$.
Note that the number of occurrences of the pattern $[1,2,3,4]$ in $\pi$ equals
the number of occurrences of that pattern in the beheaded permutation $\pi_2 \dots \pi_n$
{\bf plus} the number of the patterns $[1,2,3]$ in the beheaded permutation $\pi_2 \dots \pi_n$
where the ``1'' is $i+1$, or $i+2$, or $\dots$ or $n$.  
Let $\pi'$ be the reduction to $\{1, \dots , n-1\}$ of 
that beheaded permutation.
Also note that the number of occurrences of the pattern $[1,2,3]$ where the ``$1$'' is an $i$ gets
increased by the number of occurrences of the pattern $[1,2]$ in the beheaded permutation,
where the ``$1$'' is a $j$ with $j>i$.
We see that
$$
weight(\pi)=y_i^{n-i} weight(\pi')
\,\, 
\vert_{\, x_i \rightarrow qx_{i+1}
\,\, , \,\,  x_{i+1} \rightarrow qx_{i+2}, \,\, \dots \,\, ,  \,\, x_{n-1} \rightarrow qx_{n} \quad ; \quad
y_i \rightarrow x_iy_{i+1}
\,\, , \,\,  y_{i+1} \rightarrow x_iy_{i+2}, \,\, \dots \,\, ,  \,\, y_{n-1} \rightarrow x_iy_{n} 
}  \quad .
$$
The factor of $y_i^{n-i}$ is because converting $\pi'$ from a permutation of $\{1, \dots, n-1\}$ to
a permutation of $\{1, \dots, i-1, i+1, \dots ,n \}$, and sticking an $i$ at the front
introduces $n-i$ new $[1,2]$ patterns where the ``$1$'' is $i$.
This gives the {\it Noonan-Zeilberger Functional Equation}  for the pattern $[1,2,3,4]$:
$$
P_n(q; x_1 , \dots , x_n \, ; \, y_1 , \dots , y_n )=
\sum_{i=1}^{n} y_i^{n-i} P_{n-1}(q; x_1, \dots, x_{i-1}, qx_{i+1}, \dots, qx_n \, ; \, 
y_1, \dots, y_{i-1}, x_iy_{i+1}, \dots, x_iy_n \,
) \quad .
\eqno(NZFE2)
$$
Having found $P_n(q;x_1, \dots, x_n \, ; \, y_1, \dots, y_n )$, we set the ``catalytic'' variables $x_1, \dots, x_n$ 
and $y_1, \dots, y_n$ all to $1$
and get
$$
g_n(q):=A_{[1,2,3,4]}(q,n)=P_n(q; 1, 1, \dots, 1 \, ; \, 1, 1, \dots, 1 ) \quad.
$$
Even though this is an ``exponential-time'' (and memory!) algorithm, it is still faster than
the direct weighted counting of all the $n!$ permutations, and we were able to
explicitly compute them through $n=10$. 

The first few polynomials are
$$
g_1(q)=1 \quad, \quad  g_2(q)=2 \quad, \quad  g_3(q)=3 \quad, \quad g_4(q)=q+ 23 \quad ,
$$
$$
g_5(q)={q}^{5}+4\,{q}^{2}+12\,q+103 \quad ,
$$

$$
g_6(q)={q}^{15}+5\,{q}^{9}+8\,{q}^{6}+12\,{q}^{5}+6\,{q}^{4}+10\,{q}^{3}+63\,{q}^{2}+102\,q+513 \quad ,
$$

$$
g_7(q)={q}^{35}+6\,{q}^{25}+10\,{q}^{19}+18\,{q}^{16}+12\,{q}^{15}+13\,{q}^{13}+24\,{q}^{11}+32\,{q}^{10}+72\,{q}^{9}+10\,{q}^{8}+46\,{q}^{7}
$$
$$
+142\,{q}^{6}+116\,
{q}^{5}+146\,{q}^{4}+196\,{q}^{3}+665\,{q}^{2}+770\,q+2761 \quad ,
$$

$$
g_8(q)=
{q}^{70}+7\,{q}^{55}+12\,{q}^{45}+15\,{q}^{41}+10\,{q}^{39}+8\,{q}^{36}+28\,{q}^{35}+40\,{q}^{32}+41\,{q}^{29}+10\,{q}^{28}+24\,{q}^{27}+44\,{q}^{26}
+84\,{q}^{25}
$$
$$
+24\,{q}^{24}+89\,{q}^{23}+12\,{q}^{21}+142\,{q}^{20}+136\,{q}^{19}+96\,{q}^{18}+115\,{q}^{17}+333\,{q}^{16}+156\,{q}^{15}+112\,{q}^{14}
+312\,{q}^{13}
$$
$$
+199\,{q}^{12}+600\,{q}^{11}+573\,{q}^{10}+804\,{q}^{9}+503\,{q}^{8}+885\,{q}^{7}+1782\,{q}^{6}+1204\,{q}^{5}+2148\,{q}^{4}+2477\,{q}^{3}+5982\,{q
}^{2}+5545\,q+15767 \quad .
$$

For $g_9(q), g_{10}(q)$ see: \quad {\tt http://www.math.rutgers.edu/\~{}zeilberg/tokhniot/oP1234d}.

The obvious analog of the Crucial Lemma still holds, and one can get
{\it polynomial time} (in $n$) algorithms, to compute the number of permutations
of length $n$ with exactly $r$ occurrences of the the pattern $[1,2,3,4]$.
Alas, because we have twice as many catalytic variables, the $O(n^{r+3})$ becomes
$O(n^{2r+5})$. Nevertheless, we were able to compute the first $70$ terms for
the case $r=1$. Here are the first $23$ terms:

$$
0, 0, 0, 1, 12, 102, 770, 5545, 39220, 276144, 1948212, 13817680, 98679990, 
$$
$$
710108396, 5150076076, 37641647410, 277202062666, 2056218941678, 15358296210724,
$$
$$
115469557503753, 873561194459596, 6647760790457218, 50871527629923754 \quad .
$$

The rest can be viewed in: {\tt http://www.math.rutgers.edu/\~{}zeilberg/tokhniot/oF1234a}.

{\bf The Maple package P1234}

Everything is implemented in the Maple package {\tt P1234}. See the webpage of this article for
sample input and output files. The package {\tt F1234} is a more efficient implementation
for a small number of occurrences $r$.

{\bf Beyond}

Of course, the same reasoning applies to any increasing pattern $[1, \dots, k]$ but we have,
in addition to $q$, $(k-2)n$ additional catalytic variables. For each
specific $r$ this implies a scheme that enables one to compute  in ``polynomial'' time 
(in $n$, but of course not in $k$ or $r$) the desired numbers. For the patterns
$[1,2,3,4,5]$ and $[1,2,3,4,5,6]$ (i.e. $k=5$ and $k=6$) this is implemented in
Maple packages {\tt P12345} (and its more efficient [for small $r$] version {\tt F12345})
and {\tt P123456} respectively.

{\bf Other Patterns}

Even the case of {\it pattern-avoidance}, i.e. $r=0$, is already extremely difficult in general.
As we mentioned above for the pattern $[1,3,2,4]$, there is no known polynomial time algorithm
for enumerating permutations that avoid it. But for some few infinite families (see [Bo3] and [Wiki])
exact formulas (for the avoiding, $r=0$, case) are known, and the present approach would hopefully
be able to find polynomial-time schemes for $r>0$,
at least for some of them. We hope to investigate this in a future paper.

{\bf References}

[Bo1] Mikl\'os B\'ona,
{\it Permutations with one or two 132-subsequences},
Discrete Mathematics {\bf 175}(1997), 55-67.

[Bo2] Mikl\'os B\'ona,
{\it The Number of Permutations with Exactly
$r$ 132-Subsequences Is P-Recursive in the Size!},
Adv. Appl. Math. {\bf 18} (1997), 510-522.

[Bo3] Mikl\'os B\'ona,
{\it ``Combinatorics of Permutations''}, Chapman and Hall, 2004.

[Bo4] Mikl\'os B\'ona,
{\it The copies of any permutation pattern are asymptotically normal}, \hfill\break
{\tt http://arxiv.org/abs/0712.2792} (17 Dec. 2007).

[Bu] Alexander Burstein, {\it  A short proof for the number of permutations containing pattern 321 exactly once},
Electron. J. Combin. {\bf}18(2) (2011), \#21, (3 pp). 

[FM] Ghassan Firro and Toufik Mansour,
{\it Three-letter-pattern-avoiding permutations and
functional equations},
Electron. J. Combin. {\bf 13}(1) (2006), \#51, (14 pp).

[F] Markus Fulmek,
{\it Enumeration of permutations containing a 
prescribed number of occurrences of a pattern of length three},
Adv. Appl. Math. {\bf 30} (2003), 607-632

[K] Christoph Koutschan. 
{\it HolonomicFunctions (User's Guide)}. Technical report no. 10-01 in RISC Report Series, University of Linz, Austria.  January 2010. \hfill\break
{\tt http://www.risc.jku.at/publications/download/risc\_3934/hf.pdf}

[MV] Toufik Mansour and Alek Vainshtein,
{\it Counting occurrences of 132 in a permutation},
Adv. in Appl. Math. {\bf 28}(2002), 185-195.

[N] John Noonan,
{\it The number of permutations containing exactly one increasing 
subsequence of length three},
Discrete Mathematics {\bf 152}(1996), 307-313.

[NZ] John Noonan and Doron Zeilberger,
{\it The enumeration of permutations with a prescribed number of `forbidden' patterns},
Adv. Appl. Math. {\bf 17}(1996), 381-407.\hfill\break
http://www.math.rutgers.edu/\~{}zeilberg/mamarim/mamarimhtml/forbid.html

[SiSc] Rodica Simion and Frank W. Schmidt, {\it Restricted permutations},  European Journal of Combinatorics {\bf 6} (1985), 383-406.

[Wiki]
{\it Enumerations of specific permutation classes}, \hfill\break
{\tt http://en.wikipedia.org/wiki/Enumerations\_of\_specific\_permutation\_classes} \hfill\break 
[article initiated by Vince Vatter]

[Wil] Herbert Wilf, {\it What is an answer?}, Amer. Math. Monthly {\bf 89} (1982), 289-292.

[Z1] Doron Zeilberger,
{\it
Alexander Burstein's Lovely Combinatorial 
Proof of John Noonan's Beautiful Formula that the number 
of n-permutations that contain the Pattern 321 Exactly Once 
Equals (3/n)(2n)!/((n-3)!(n+3)!)},
Personal Journal of Shalosh B. Ekhad and Doron Zeilberger,
Oct. 18, 2011.
\hfill\break
http://www.math.rutgers.edu/\~{}zeilberg/mamarim/mamarimhtml/burstein.html

[Z2] Doron Zeilberger, {\it A Holonomic Systems Approach To Special Functions},
J. Computational and Applied Math {\bf 32}(1990), 321-368. 
\hfill\break
http://www.math.rutgers.edu/\~{}zeilberg/mamarim/mamarimhtml/holonomic.html

\end